\documentclass[a4paper,12pt]{article}

\usepackage{amsmath}
\usepackage{amsfonts}
\usepackage{amssymb}
\usepackage{latexsym}
\usepackage{epsfig}
\usepackage{graphicx}
\usepackage{oldgerm}
\usepackage{theorem}
\usepackage{bm}

\def\C{\mathbb{C}}
\def\R{\mathbb{R}}
\def\N{\mathbb{N}}
\def\B{\mathbb{B}}
\def\Z{\mathbb{Z}}
\def\D{\mathbb{D}}

\def\bE{\mathbf{E}}

\def\bH{\mathbf{H}}

\def\cB{\mathcal{B}}

\def\cD{\mathcal{D}}
\def\cF{\mathcal{F}}

\def\sH{\mathsf{H}}
\def\sh{\mathsf{h}}
\def\sB{\mathsf{B}}

\def\Conf{\mathrm{Conf}}
\def\var{\mathrm{var}}
\def\vol{\mathrm{vol}}
\def\Re{\mathrm{Re}\,}
\def\Im{\mathrm{Im}\,}
\def\rR{\mathrm{R}}
\def\rI{\mathrm{I}}
\def\rN{\mathrm{N}}

\def\bra{\langle}
\def\ket{\rangle}
\def\dis={\stackrel{\rm d}{=}}
\def\law={\stackrel{\rm (law)}{=}}

\theorembodyfont{\itshape}
\newtheorem{thm}{Theorem}[section]
\newtheorem{lem}[thm]{Lemma}

\newtheorem{prop}[thm]{Proposition}

\newtheorem{df}[thm]{Definition}

\newcommand{\SSC}[1]{\section{#1}\setcounter{equation}{0}}

\begin{document}

\title{
Hyperuniformity of the 
determinantal point processes \\
associated with the Heisenberg group
\footnote{
This manuscript was prepared for the Proceedings 
of the 2021 RIMS (Research Institute for
Mathematical Sciences) Workshop 
`Mathematical Aspects of Quantum Fields
and Related Topics', which was held online
on December 6--8, 2021. 
The Proceedings will be issued in RIMS K\^{o}ky\^{u}roku
edited by Fumio Hiroshima
and published at Kyoto University Research Information
Repository and RIMS Homepage.
}
}
\author{
Makoto Katori 
\footnote{
Department of Physics,
Faculty of Science and Engineering,
Chuo University, 
Kasuga, Bunkyo-ku, Tokyo 112-8551, Japan;
e-mail: katori@phys.chuo-u.ac.jp
} 
}

\date{17 March 2022}
\pagestyle{plain}
\maketitle

\begin{abstract}
The Ginibre point process is given by the eigenvalue
distribution of a non-hermitian complex Gaussian matrix
in the infinite matrix-size limit.
This is a determinantal point process (DPP) on the complex
plane ${\mathbb{C}}$ in the sense that all correlation functions
are given by determinants specified by an integral
kernel called the correlation kernel.
Shirai introduced the one-parameter
($m \in {\mathbb{N}}_0$) extensions of the Ginibre DPP 
and called them
the Ginibre-type point processes.
In the present paper we consider a generalization
of the Ginibre and the Ginibre-type point processes 
on ${\mathbb{C}}$
to the DPPs in the higher-dimensional spaces, 
${\mathbb{C}}^D, D=2,3, \dots$, in which
they are parameterized by a multivariate level
$m \in {\mathbb{N}}_0^D$. 
We call the obtained point processes 
the extended Heisenberg family of DPPs,
since the correlation kernels are generally 
identified with the correlations of two points in
the space of Heisenberg group expressed by 
the Schr\"{o}dinger representations.
We prove that all DPPs in this large family are 
in Class I of hyperuniformity. 

\vskip 0.5cm


\end{abstract}
\vspace{3mm}
\normalsize

\SSC
{Introduction and Results} \label{sec:Introduction}

We consider the $d$-dimensional Euclid space
$\R^d$, $d \in \N:=\{1,2, \dots\}$, 
or the $D$-dimensional complex space $\C^D$, $D \in \N$
as a base space $S$. 
We assume that $S$ is associated with a reference measure
$\lambda$.
We consider an \textit{infinite point process} on $S$,
which is expressed by an infinite sum of delta measures 
concentrated on a set of random points
$X_i, i \in \N$,
\begin{equation}
\Xi = \sum_{i: i \in \N} \delta_{X_i}.
\label{eqn:Xi}
\end{equation}
We assume that for any 
bounded domain $\Lambda \subset S$,
$\Xi(\Lambda) < \infty$; 
that is, accumulation of points does not occur.
We also assume that 
with respect to the reference measure $\lambda(dx)$
the point process has \textit{a finite density
$\rho_1(x) < \infty$ }
at almost every $x \in S$.

We consider a \textit{homogeneous} point process
in the sense that 
\[
\rho_1(x) \lambda(dx)= {\rm const.} \times dx, \quad
x \in S, 
\]
where $dx$ denotes the Lebesgue measure on $S$.
The above assumption implies that 
for a bounded domain $\Lambda \subset S$, 
\[
\bE[\Xi(\Lambda)] \propto \vol(\Lambda).
\]
Now we consider the \textit{number variance}
in the domain $\Lambda$,
\[
\var[\Xi(\Lambda)] := \bE[(\Xi(\Lambda)- \bE[\Xi(\Lambda)])^2],
\]
which represents local density fluctuation of 
the point process $\Xi$. 
The domain $\Lambda$ is regarded as
an \textit{observation window} to measure 
the density fluctuation.
If the points are non-correlated and
given by a Poisson point process, then 
\[
\var[\Xi(\Lambda)] \propto \vol(\Lambda).
\]

Recently in condensed matter physics and related material sciences,
correlated particle systems are said to be in a
\textit{hyperuniform state} when density fluctuations
are anomalously suppressed in large-scale limit.
(See, for instance, \cite{GL17,Tor18}.) 
For an infinite random point process $\Xi$, 
the \textit{hyperuniformity} is defined by 
\[
\lim_{\Lambda \to S}
\frac{\var[\Xi(\Lambda)]}{\bE[\Xi(\Lambda)]}
=0. 
\]
This means that 
the number variance of points grows
more slowly than the window volume in the limit
such that the window covers whole of the space
$\Lambda \to S$. 

Torquato \cite{Tor18} proposed
\textit{three hyperuniformity classes} for point processes
concerning asymptotics of number variances.
In order to clearly assert this classification,
here we assume that $S=\R^d$, $d \in \N$, and
$\Lambda=\B^{(d)}_R:=\{x \in \R^d : |x| < R\}$, $R>0$, 
where
$\vol(\B^{(d)}_R) = \pi^{d/2} R^d/\Gamma(d/2+1)$
with the gamma function
$\Gamma(z) :=\int_0^{\infty} e^{-u} u^{z-1} du$, 
$\Re z >0$.
We consider a series of balls with increasing radius $R$,
$\{\B^{(d)}_R \}_{R >0}$, and
the hyperuniform states are classified 
as follows; 
\begin{align*}
&\mbox{Class I} : \qquad \, \, \,
\var[\Xi(\B^{(d)}_R)]
\asymp R^{d-1},
\nonumber\\
&\mbox{Class II} : \qquad \, \,
\var[\Xi(\B^{(d)}_R)]
\asymp R^{d-1} \log R,
\nonumber\\
&\mbox{Class III} : \qquad
\var[\Xi(\B^{(d)}_R)]
\asymp R^{d-\alpha}, \quad 0<\alpha< 1,
\quad \mbox{as $R \to \infty$}.
\end{align*}
Here $f(R) \asymp g(R)$ means that there
are finite positive constants $c_1$ and $c_2$ such that
$c_1 g(R) < f(R) < c_2 g(R)$.
The above characterization of three classes
will be similarly described for any 
series of windows $\{\Lambda_R \}_{R >0}$
labeled by a linear scale $R$ of window.
It is expected that the hyperuniformity and its classification
are the proper properties of $\Xi$ and do not depend on
the choice of observation window $\Lambda_R$, $R >0$.

Determinantal point processes (DPPs)
studied in random matrix theory (RMT) \cite{For10}
provide a variety of examples of hyperuniform systems.
In general a DPP is specified by a triplet 
\[
(\Xi, K, \lambda(dx)), 
\]
where
$\Xi$ is a nonnegative-integer-valued
Radon measure (\ref{eqn:Xi}) 
representing a point process, 
$K$ is a continuous function
$S \times S \to \C$ called the
\textit{correlation kernel},
and $\lambda(dx)$ is a reference measure on $S$. 

The most studied DPP in RMT may be the
\textit{sinc} (\textit{sine}) \textit{DPP}, 
$(\Xi_{\rm sinc}, K_{\rm sinc}, dx)$ on $S=\R$
with the correlation kernel
\[
K_{\rm sinc}(x,y)=
\frac{\sin(x-y)}{\pi(x-y)}, \quad x, y \in \R. 
\]
This DPP is obtained as the bulk scaling limit
of the eigenvalue distribution of \textit{Hermitian random matrices}
in the Gaussian unitary ensemble (GUE) \cite{For10}.
If $\lim_{R \to \infty} f(R)/g(R)=1$, we will write
$f(R) \sim g(R)$ as $R \to \infty$. 
As a classical result in RMT, 
it is well known that
\[
\var[\Xi_{\rm sinc}(\B^{(1)}_R)] 
\sim 
\frac{\log R}{\pi^2} 
\quad \mbox{as $R \to \infty$}.
\]
(See, for instance, \cite{CL95,Sos00b,Sos02}, 
\cite[Remark 5.8]{ST03}.)
That is, the sinc DPP is in Class II of hyperuniformity.
Torquato \cite{Tor18}
studied one-parameter ($d \in \N$) family of DPPs
called the \textit{Fermi-sphere point processes},
which is also called 
the \textit{Euclidean family of DPPs} in \cite{KS21}.
This family gives the sinc DPP when $d=1$.
It was proved that the Fermi-sphere point processes
are in Class II of hyperuniformity for general 
$d \in \N$ \cite{TSZ08,Tor18}.

An example of infinite DPP in 
Class I of hyperuniformity is also provided in RMT \cite{For10}.
It is the DPP on $\C$ called the
\textit{Ginibre DPP}, 
$(\Xi_{\rm Ginibre}, K_{\rm Ginibre}, \lambda_{\rN(0,1; \C)}(dx))$
on $S=\C$, 
which is obtained as the
bulk scaling limit of eigenvalue distribution of 
\textit{non-Hermitian random matrices} \cite{Gin65}. There
\begin{align*}
K_{\rm Ginibre}(x,y) &=e^{x \overline{y}},
\quad x, y \in \C, \nonumber\\
\lambda_{\rN(0,1;\C)}(dx) &= \frac{e^{-|x|^2}}{\pi} dx.
\end{align*}
Note that the disk on $\C$, 
$\D_R :=\{x \in \C: |x| < R \}$, is identified with 
$\B^{(2)}_R \subset \R^2$.
Shirai proved \cite{Shi06}
\begin{equation}
\var[\Xi_{\rm Ginibre}(\B^{(2)}_R) ]
\sim \frac{R}{\sqrt{\pi}} \quad
\mbox{as $R \to \infty$}.
\label{eqn:Shirai06}
\end{equation}

In \cite{MKS21} 
the result (\ref{eqn:Shirai06}) 
in $S=\C$ was extended to DPPs in the
higher-dimensional complex spaces
$S=\C^D, D=2,3 \dots$ as follows.
For $S=\C^D, D \in \N$, each coordinate
$x \in \C^D$ has $D$ complex components;
$x=(x^{(1)}, \dots, x^{(D)})$ with
$x^{(\ell)}=\Re x^{(\ell)}+ \sqrt{-1} \Im x^{(\ell)}$, $\ell=1, \dots, D$.
We set $x_{\rR}:=(\Re x^{(1)}, \dots, \Re x^{(D)})$,
$x_{\rI}:=(\Im x^{(1)}, \dots, \Im x^{(D)}) \in \R^D$,
and we write $x=x_{\rR}+ \sqrt{-1} x_{\rI}$ in this paper.
The Lebesgue measure on $\C^D$ is given by
$d x = d x_{\rR} d x_{\rI} 
:= \prod_{\ell=1}^D d \Re x^{(\ell)} d \Im x^{(\ell)}$.
For $x = x_{\rR}+ \sqrt{-1} x_{\rI}$,
$y=y_{\rR}+ \sqrt{-1} y_{\rI} \in \C^D$,
we use the \textit{standard Hermitian inner product}; 
\begin{align*}
x \cdot \overline{y} 
&:= (x_{\rR}+ \sqrt{-1} x_{\rI}) \cdot (y_{\rR}-\sqrt{-1} y_{\rI})
\nonumber\\
&=(x_{\rR} \cdot y_{\rR}+x_{\rI} \cdot y_{\rI})
-\sqrt{-1} (x_{\rR} \cdot y_{\rI}- x_{\rI} \cdot y_{\rR}).
\end{align*}
The norm is given by
$|x| := \sqrt{x \cdot \overline{x}} = \sqrt{|x_{\rR}|^2+|x_{\rI}|^2}$,
$x \in \C^D$.
Notice that if $x=x_{\rR}, y=y_{\rR} \in \R^D$,
then 
$x \cdot \overline{y}=x_{\rR} \cdot y_{\rR} 
:=\sum_{\ell=1}^D \Re x^{(\ell)} \Re y^{(\ell)}$. 
The $D$-dimensional disk 
$\D^{(D)}_R:=\{x \in \C^D: |x| < R\}$ is
identified with $\B^{(d)}_{R}$ in $\R^d$ 
provided that $d=2D, D \in \N$.
On $\C^D$ the reference measure is given by
the $D$-dimensional direct-product of $\lambda_{\rN(0,1; \C)}(dx)$, 
\begin{align}
\lambda_{\rN(0, 1; \C^D)}(d x)
&:= \prod_{\ell=1}^D \lambda_{\rN(0, 1; \C)}(d x^{(\ell)})
\nonumber\\
&= \frac{e^{-(|x_{\rR}|^2+|x_{\rI}|^2)}}{\pi^D} dx_{\rR} dx_{\rI}
= \frac{e^{-|x|^2}}{\pi^D} dx, 
\quad x \in \C^D.
\label{eqn:lambda1}
\end{align}
In \cite{MKS21} 
the one-parameter family ($D \in \N$) of DPPs
was studied, which is called the
\textit{Heisenberg family of DPPs} 
defined on $\C^D$ as follows.

\begin{df}
\label{thm:HeisenbergDPP}
The Heisenberg family of DPPs 
is defined by
$(\Xi_{\sH_D}, K_{\sH_D}, \lambda_{\rN(0, 1; \C^D)}(dx))$
on $\C^D$, $D \in \N$
with the correlation kernel 
\begin{equation}
K_{\sH_D}(x, y)
=e^{x \cdot \overline{y}},
\quad x, y \in \C^D.
\label{eqn:K_HD}
\end{equation}
\end{df}

Note that $K_{\sH_D}$ is hermitian;
$\overline{K_{\sH_D}(x,y)}=K_{\sH_D}(y,x)$, $x, y \in \C^D$. 
The kernels in this form on $\C^D, D \in \N$ have been
studied by Zelditch,
who identified them with the Szeg\H{o} kernels
for the \textit{reduced Heisenberg group} $\sH_D^{\rm red}$
\cite{Zel01}.
This family includes the Ginibre DPP as the
lowest dimensional case; $D=1$.
The following was proved \cite{MKS21}.

\begin{thm}
\label{thm:main2}
\textit{Any DPP in the Heisenberg family}, 
$(\Xi_{\sH_D}, K_{\sH_D}, \lambda_{\rN(0, 1; \C^D)}(dx))$
on $\C^D$, $D \in \N$, 
is in Class I of hyperuniformity such that
\[
\lim_{R \to \infty} R \frac{\var[\Xi_{\sH_D}(\B^{(2D)}_R)]}
{\bE[\Xi_{\sH_D}(\B^{(2D)}_R)]}
=\frac{D}{\sqrt{\pi}}.
\label{eqn:C0}
\]
Moreover, for each $D \in \N$, 
the following \textit{asymptotic expansion} holds, 
\[
\frac{\var[\Xi_{\sH_D}(\B^{(2D)}_R)]}
{\bE[\Xi_{\sH_D}(\B^{(2D)}_R)]}
\sim
\frac{D}{\sqrt{\pi}} R^{-1}
\sum_{k=0}^{\infty} (-1)^k \frac{\alpha_k(D)}{(2k+1) k! 2^{4k}} R^{-2k}
\quad \mbox{as $R \to \infty$}, 
\label{eqn:asym}
\]
where 
\[
\alpha_k(D) = 
\begin{cases}
1, & \quad \mbox{if $k=0$},
\cr
\displaystyle{
\prod_{\ell=-k+1}^k(2D+2 \ell-1)
},
& \quad \mbox{if $k \in \N$}. 
\end{cases}
\]
\end{thm}

Instead of $\B^{(2D)}_R \simeq \D^{(D)}_R$, $D \in \N$,
we can consider the \textit{$D$-dimensional polydisk}
of radius $R>0$ in $\C^D$, 
\begin{equation}
\Delta^{(D)}_R:= \{ x=(x^{(1)}, \dots, x^{(D)}) \in \C^D :
|x^{(i)}| < R, i=1, \dots, D \},
\label{eqn:polydisk}
\end{equation}
as an observation window. 
In Remark 5 of \cite{MKS21}, the following was proved.

\begin{prop}
\label{thm:MKS_R5}
For the Heisenberg family of DPPs, 
$(\Xi_{\sH_D}, K_{\sH_D}, \lambda_{\rN(0,1; \C^D)}(dx))$,
$D \in \N$, 
\begin{align*}
\lim_{R \to \infty} R \frac{\var[\Xi_{\sH_D}(\Delta^{(D)}_R)]}
{\bE[\Xi_{\sH_D}(\Delta^{(D)}_R)]}
&=\lim_{R \to \infty} R 
\left[
1- \left(
1-\frac{\var[\Xi_{\sH_1}(\B^{(2)}_R)]}
{\bE[\Xi_{\sH_1}(\B^{(2)}_R)]} \right)^D
\right]
\nonumber\\
&= \frac{D}{\sqrt{\pi}}. 
\end{align*}
\end{prop}

The 
\textit{Laguerre polynomial} is defined by
\begin{equation}
L_n^{(\alpha)}(\zeta) :=
\frac{\zeta^{-\alpha} e^{\zeta}}{n!}
\frac{d^n}{d \zeta^n} (\zeta^{n+\alpha} e^{-\zeta}),
\quad n \in \N_0:=\{0, 1, \dots\}, \quad \alpha, 
\zeta \in \R.
\label{eqn:Laguerre}
\end{equation}
In particuler, we wite
$L_n(\zeta) := L_n^{(0)}(\zeta)$.
Note that $L_0^{(\alpha)}(\zeta)=1$ for
any $\alpha \in \R$.
Let $_3F_2(\alpha_1, \alpha_2, \alpha_3; \beta_1, \beta_2; z)$
be the hypergeometric function defined by 
\[
_3F_2(\alpha_1, \alpha_2, \alpha_3; \beta_1, \beta_2; z)
:=\sum_{n=0}^{\infty}
\frac{(\alpha_1)_n (\alpha_2)_n (\alpha_3)_n}
{(\beta_1)_n (\beta_2)_n} \frac{z^n}{n!}, 
\]
where $(\alpha)_n$ is the Pochhammer symbol; 
$(\alpha)_n:=\alpha(\alpha+1) \cdots (\alpha+n-1), n \in \N$,
$(\alpha)_0 := 1$.
It is obvious that
\begin{equation}
_{3}F_2(\alpha_1, \alpha_2, 0; \beta_1, \beta_2; z)=1.
\label{eqn:hypergeo1}
\end{equation}
For $D=1$, let
\[
K_{\sH_1}^{(m)}(x, y) 
=K_{\sH_1}(x, y) L_m(|x-y|^2),
\quad m \in \N_0, \quad x, y \in \C.
\]
The one-parameter ($m \in \N_0$) family of DPPs, 
$(\Xi_{\sH_1}^{(m)}, K_{\sH_1}^{(m)}, \lambda_{\rN(0,1; \C)})$,
$m \in \N_0$ on $\C$ was studied by Shirai \cite{Shi15},
who called the DPPs in this family the
\textit{Ginibre-type point processes}.
This family of DPPs was also studied by
Haimi and Hedenmalm \cite{HH13}, where they called the DPPs
the \textit{polyanalytic Ginibre ensembles}.
This family of DPPs on $\C$ has 
a physical interpretation in terms of
a two-dimensional system of free electrons 
in a uniform magnetic field,
where the electrons occupy the
\textit{$m$-th Landau energy level},
$m \in \N_0$. 
Shirai proved the following \cite{Shi15}.

\begin{prop}
\label{thm:Shirai1}
Any \textit{Ginibre-type DPP on $\C$},
$(\Xi_{\sH_1}^{(m)}, K_{\sH_1}^{(m)}, \lambda_{\rN(0,1; \C)}(dx))$,
$m \in \N_0$, 
is in Class I of hyperuniformity
such that
\[
\lim_{R \to \infty} R \frac{\var[\Xi_{\sH_1}^{(m)}(\B^{(2)}_R)]}
{\bE[\Xi_{\sH_1}^{(m)}(\B^{(2)}_R)]}
=C_{\sH_1}^{(m)}
\label{eqn:Shirai1}
\]
with 
\begin{align}
C_{\sH_1}^{(m)}
&=\frac{2 \Gamma(m+3/2)}{\pi m!} 
{_{3}F_2} 
\left(- \frac{1}{2}, -\frac{1}{2}, -m; 1, - \frac{1}{2}-m; 1 \right)
\nonumber\\
&\sim \frac{8}{\pi^2} m^{1/2} 
\quad \mbox{as $m \to \infty$}.
\label{eqn:CH1}
\end{align}
\end{prop}

In the present paper, we consider 
a large family of DPPs which includes both of
the Heisenberg family of DPPs and the Ginibre-type DPPs, 
which is parameterized by the dimensionality $D$
of the base space $S=\C^D$ and the \textit{multivariate level}
expressed by a set of integers,
$m=(m^{(1)}, \dots, m^{(D)}) \in \N_0^D$. 

\begin{df}
\label{thm:extended_HeisenbergDPP}
The extended Heisenberg family of DPPs 
is defined by
$(\Xi^{(m)}_{\sH_D}$, $K_{\sH^{(m)}_D}$, 
$\lambda_{\rN(0, 1; \C^D)}(dx))$
on $\C^D$, $D \in \N$
with a multivariate level $m \in \N_0^D$, where
the correlation kernel is given by
\begin{equation}
K^{(m)}_{\sH_D}(x, y)
=K_{\sH_D}(x, y) \prod_{\ell=1}^D 
L_{m^{(\ell)}} (|x^{(\ell)}-y^{(\ell)}|^2), 
\quad x, y \in \C^D.
\label{eqn:K_HD_m}
\end{equation}
Here $K_{\sH_D}$ is defined by (\ref{eqn:K_HD}). 
\end{df}

When $m=0:=(0, \dots, 0)$, this family is
reduced to the original Heisenberg family
defined by Theorem \ref{thm:HeisenbergDPP}.
It is obvious that when $D=1$ this family is
identified with the Ginibre-type DPPs of Shirai \cite{Shi15}.

In \cite{Shi15}
Shirai gave a sufficient condition of Class I of hyperuniformity
for DPPs on $\R^d$ with general dimensions $d \geq 1$.
There he assumed that the observation windows are balls
$\B^{(d)}_R $
and the DPPs are \textit{isotropic} in a sense.
We can apply this general theorem to
the Heisenberg family of DPPs,  
$(\Xi_{\sH_D}, K_{\sH_D}, \lambda_{\rN(0,1; \C^D)}(dx))$
for all $D \in \N$ and
the Ginibre-type DPPs on $\C$, 
$(\Xi_{\sH_1}^{(m)}, K_{\sH_1}^{(m)}, \lambda_{\rN(0,1; \C)}(dx))$
for all $m \in \N_0$.
But, it is not applicable to
the extended Heisenberg family of DPPs on $\C^D$,
$(\Xi_{\sH_D}^{(m)}, K_{\sH_D}^{(m)}, \lambda_{\rN(0,1; \C^D)}(dx))$, 
if $D \geq 2$ and $m \not=0$, since
\[
\left|
\sqrt{ \frac{e^{-|x|^2}}{\pi^D}}
K_{\sH_D}^{(m)}(x, x') 
\sqrt{ \frac{e^{-|x'|^2}}{\pi^D}} \right|^2
=\frac{1}{\pi^{2D}} e^{-|x-x'|^2}
\prod_{\ell=1}^D L_{m^{(\ell)}}(|x^{(\ell)}-{x'}^{(\ell)}|^2)^2
\]
is not a function of 
$|x-x'|^2 =\sum_{\ell=1}^n |x^{(\ell)}-{x'}^{(\ell)}|^2$.
In other words, the DPPs in the extended Heisenberg family
with $D \geq 2$ or $m \not=0$ are not isotropic.

In the present paper, however, 
we prove the following result
for the extended Heisenberg family of DPPs,
when observation windows are given by
the $D$-dimensional polydisks (\ref{eqn:polydisk}).

\begin{thm}
\label{thm:polydisk1}
For the extended Heisenberg family of DPPs, 
$(\Xi_{\sH_D}^{(m)}, K_{\sH_D}^{(m)}, \lambda_{\rN(0,1; \C^D)}(dx))$,
$D \in \N$, 
$m \in \N_0^D$, 
\begin{align}
\lim_{R \to \infty} R \frac{\var[\Xi_{\sH_D}^{(m)}(\Delta^{(D)}_R)]}
{\bE[\Xi_{\sH_D}^{(m)}(\Delta^{(D)}_R)]}
&=\lim_{R \to \infty} R 
\left[
1- \prod_{\ell=1}^D \left(
1-\frac{\var[\Xi_{\sH_1}^{(m^{(\ell)})}(\B^{(2)}_R)]}
{\bE[\Xi_{\sH_1}^{(m^{(\ell)})}(\B^{(2)}_R)]} \right)
\right]
\nonumber\\
&= \sum_{\ell=1}^D C_{\sH_1}^{(m^{(\ell)})}.
\label{eqn:Main1}
\end{align}
\end{thm}
By (\ref{eqn:hypergeo1}), (\ref{eqn:CH1}) gives
\[
C_{\sH_1}^{(0)}= \frac{2 \Gamma(3/2)}{\pi}
=\frac{1}{\sqrt{\pi}}.
\]
Hence, this result can be regarded as a generalization
of Proposition \ref{thm:MKS_R5}.

In Section \ref{sec:Heisenberg}, we show 
the relationship between 
the extended Heisenberg family of DPPs and 
the \textit{Schr\"{o}dinger representations of
the Heisenberg group $\sH_D$},
$D \in \N$. 
There the function denoted by $g$
specifies the representation.
The correlation kernels 
$K_{\sH_D}^{(m)}$, $D \in \N$, $m \in \N_0^D$ are
identified with the correlations of two points
in the space of $\sH_D$ expressed using the 
Schr\"{o}dinger operators acting on $g$.
In Section \ref{sec:preliminaries_proofs}, 
after giving preliminaries for point processes
and comments on 
Theorem \ref{thm:main2}, 
the proof of Theorem \ref{thm:polydisk1} is given.

By Abreu and his coworkers \cite{AGR16,AGR19}, 
the DPPs associated with the 
Schr\"{o}dinger representations of
$\sH_D$ are named as
the \textit{Weyl--Heisenberg ensembles},
in which $g$ is called 
a \textit{window function} in the context of the 
\textit{time-frequency analysis} \cite{Gro01}.
See also Section 2.6 of \cite{KS21}.
It was proved by Abreu et al.~\cite{APRT17}
that the Weyl--Heisenberg ensembles
are in the hyperuniform state of Class I
for a general class of window functions.

\SSC
{Representations of the Heisenberg Group
and Correlation Kernels} 
\label{sec:Heisenberg}
\subsection{Schr\"{o}dinger representations}
\label{sec:Schrodinger}

We briefly review
the representation theory of the \textit{Heisenberg group}
\cite{Fol89,Ste93,Gro01}
in order to explain the reason why we call the
DPPs defined by Definition \ref{thm:HeisenbergDPP}
the \textit{Heisenberg family} of DPPs
and why we think that the DPPs defined by Definition
\ref{thm:extended_HeisenbergDPP}
form its extended family. 

Consider the classical and quantum kinetics of a single particle
moving in $\R^D, D \in \N$.
We note that, if $D=3k, k \in \N$, this 
represents a $k$-particle system in 
the three dimensional Euclidean space.
The \textit{phase space} is given by $\R^{2D}$ with
coordinates 
\[
(p, q)=(p_1, \dots, p_D, q_1, \dots, q_D).
\]

In order to describe the 
\textit{Heisenberg Lie algebra $\sh_D$}, we consider
$\R^{2D+1}$ with coordinates
$(p, q, \tau)=(p_1, \dots, p_D, q_1, \dots, q_D, \tau)$, 
in which a Lie bracket is given by
\[
[(p,q,\tau), (p', q', \tau')]
=(0, 0, p \cdot q'-q  \cdot p')
=(0, 0, [(p,q), (p', q')]).
\]
The symplectic form of the \textit{Lie bracket} 
$[(p,q), (p',q')]=p \cdot q'-q \cdot p'$ 
comes from
the Poisson bracket in the classical mechanics
and the commutator $[A, B] := AB-BA$ in quantum mechanics.
The \textit{Heisenberg group $\sH_D$} is the Lie group
on $\R^{2D+1}$ satisfying the group law
\[
Z Z'=Z+Z'+\frac{1}{2}[Z, Z'], \quad Z, Z' \in \R^{2D+1};
\]
that is,
\[
(p, q, \tau) (p', q', \tau')
=\left(p+p', q+q', \tau+\tau'+ \frac{1}{2}(p \cdot q'-q \cdot p') \right).
\]

Let $L^2(\R^D)$ be the set of square integrable functions on
$\R^D$, where the inner product is given by
\[
\bra f, g \ket_{L^2(\R^D)}
:= \int_{\R^D} f(\zeta) \overline{g(\zeta)} d \zeta,
\quad f, g \in L^2(\R^D)
\]
with the norm
$\|f \|_{L^2(\R^D)}:=\sqrt{\bra f, f \ket_{L^2(\R^D)}}$,
$f \in L^2(\R^D)$, 
where 
$\zeta=(\zeta^{(1)}, \dots, \zeta^{(D)}) \in \R^D$ and
$d \zeta$ denotes the Lebesgue measure on $\R^D$.
For a smooth function $f$,
we introduce operators $X^{(\ell)}$ and $\cD^{(\ell)}$ defined by
\[
(X^{(\ell)} f)(\zeta)=\zeta^{(\ell)} f(\zeta), \quad
(\cD^{(\ell)} f)(\zeta)
=\frac{1}{2 \sqrt{-1}} \frac{\partial f}{\partial \zeta^{(\ell)}} (\zeta),
\quad \ell=1, \dots, D.
\]
They satisfy the commutation relations \,
\[ 
[X^{(\ell)}, \cD^{(\ell')}]=\frac{\sqrt{-1}}{2} \delta_{\ell \ell'}, \quad
\ell, \ell' = 1, \dots, D.
\]
Note that the above will represent the
\textit{canonical commutation relations} in quantum mechanics,
$[Q^{(\ell)}, P^{(\ell')}]=\sqrt{-1} \hbar \delta_{\ell \ell'}$.
Here we should claim that the value of the Planck constant
$\hbar$ is specially chosen to be 1/2.
This corresponds to the choice of 
the reference measure on $\C^D$ as (\ref{eqn:lambda1}).
We consider a map from $\sH_D$ to the group of
unitary operators acting on $L^2(\R^D)$ defined by
\[
\rho(p, q, \tau)=e^{2 \sqrt{-1} (p \cdot \cD + q \cdot X+ \tau I)},
\]
where $\cD:=(\cD^{(1)}, \dots, \cD^{(D)})$,
$X:=(X^{(1)}, \dots, X^{(D)})$ and $I$ denotes the identity operator.
By the Baker--Campbell--Hausdorff formula, we can show that 
\[
\rho(p, q, \tau) f(\zeta)
=e^{2 \sqrt{-1} (\tau + q \cdot \zeta + p \cdot q/2)} 
f(\zeta+p),
\quad f \in L^2(\R^D).
\]
The map $\rho$ is called the \textit{Schr\"{o}dinger representation}
of $\sH_D$. 
The kernel of $\rho$ is $\{(0, 0, k \pi) : k \in \Z \} $,
since $e^{2 \pi k \sqrt{-1}}=1, k \in \Z$.
The \textit{reduced Heisenberg group}
$\sH_D^{\rm red}$ is defined by
$\sH_D^{\rm red}:=\sH_D /\{(0, 0, k \pi): k \in \Z\}$.
We consider the \textit{matrix-element functions} 
of $\rho(p, q, \tau)$ 
at $(f, g) \in L^2(\R^D)^2$, 
\begin{align}
M_{f, g}(p,q,\tau)
&:= \bra \rho(p, q, \tau) f, g \ket_{L^2(\R^D)}
= \bra \rho(-p, q, \tau) \overline{g}, \overline{f} \ket_{L^2(\R^D)}
\nonumber\\
&=e^{2 \sqrt{-1} \tau} 
\int_{\R^D} e^{2 \sqrt{-1} q \cdot \zeta}
f \left( \zeta + \frac{p}{2} \right)
\overline{g \left(\zeta-\frac{p}{2} \right)} d \zeta,
\quad f, g \in L^2(\R^2), 
\label{eqn:matrix_element}
\end{align}
which is also called the \textit{Fourier--Wigner transform}
(see, for instance, \cite[Section 1.4]{Fol89}).

\subsection{Ground-state representation}
\label{sec:ground}

We put
\begin{equation}
g(\zeta)=\left( \frac{2}{\pi} \right)^{D/4}
\frac{e^{-\zeta^2}}{\pi^{D/2}} =:G_0(\zeta),
\quad \zeta \in \R^D, 
\label{eqn:G0}
\end{equation}
and define the complex variables
\[
x=(x^{(1)}, \dots, x^{(D)})
:= p + \sqrt{-1} q
=(p^{(1)}+\sqrt{-1} q^{(1)}, \dots, p^{(D)}+\sqrt{-1} q^{(D)}) 
\in \C^D.
\]
Then (\ref{eqn:matrix_element}) becomes
\, 
\[M_{f, G_0}(p, q, \tau)
=e^{2 \sqrt{-1} \tau} \sB[f](x)
\frac{e^{-|x|^2/2}}{\pi^{D/2}}
\]
with
\begin{equation}
\sB[f](x)
:= \left( \frac{2}{\pi} \right)^{D/4}
\int_{\R^D} f(\zeta) e^{2 \zeta \cdot x-\zeta^2-x^2/2} d \zeta,
\quad f \in L^2(\R^D).
\label{eqn:Bargmann}
\end{equation}
This integral
is called the \textit{Bargmann transform}.
For $f \in L^2(\R^D)$, this integral 
converges uniformly for $x$ in any compact subset of $\C^D$,
and hence $\sB[f]$ is an entire function on $\C^D$.
We can prove that for $f_1, f_2 \in L^2(\R^D)$
\begin{equation}
\bra M_{f_1, G_0}, 
M_{f_2, G_0} \ket_{L^2(\R^{2D})}
= \bra f_1, f_2 \ket_{L^2(\R^D)}
= \bra \sB[f_1], \sB[f_2]
\ket_{L^2(\C^D, \lambda_{\rN(0, 1; \C^D)})},
\label{eqn:isometry}
\end{equation}
where $\lambda_{\rN(0,1; \C^D)}$ is given by (\ref{eqn:lambda1}) 
and
\[
\bra F_1, F_2 \ket_{L^2(\C^D, \lambda_{\rN(0, 1; \C^D)})}
:= \int_{\C^D} F_1(x) \overline{F_2(x)} 
\lambda_{\rN(0,1; \C^D)}(dx)
\]
with the norm
$\| F \|_{L^2(\C^D, \lambda_{\rN(0, 1; \C^D)})}
:=\sqrt{\bra F, F \ket_{L^2(\C^D, \lambda_{\rN(0, 1; \C^D)})}}$.
The Bargmann--Fock space $\cF_D$ is
defined by
\[
\cF_D := \left\{F : \mbox{$F$ is entire on $\C^D$ and
$\| F \|_{L^2(\C^D, \lambda_{\rN(0, 1; \C^D)})} < \infty$} \right\}.
\]
The equalities (\ref{eqn:isometry}) imply that
the Bargmann transform is an 
\textit{isometry} from $L^2(\R^D)$ into $\cF_D$.
Hence, 
if $\{f_n\}_{n \in \N_0}$ is a complete orthonormal system (CONS) 
of $L^2(\R^D)$,
then $\{\sB[f_n]\}_{n \in \N_0^D}$ makes a CONS for $\cF_D$.

For $(p,q,\tau), (p', q', \tau') \in \R^{2D+1}$, 
we consider the correlation of these two points 
expressed by the Schr\"{o}dinger operators 
acting on $G_0$ such that
\begin{equation}
{\rm Cor}_{G_0}((p,q,\tau), (p', q', \tau'))
:=\bra \rho(-p, q, \tau) G_0, \rho(-p', q', \tau') G_0 \ket_{L^2(\R^D)}.
\label{eqn:CorG01}
\end{equation}
We take a CONS $\{f_n\}_{n \in \N_0}$ of $L^2(\R^D)$,
where we assume $f_n \in \R$,
$n \in \N_0$. Then
the two-point correlation (\ref{eqn:CorG01}) 
is expanded and expressed using the matrix-element functions
(\ref{eqn:matrix_element})
as
\begin{align}
{\rm Cor}_{G_0}((p,q,\tau), (p', q', \tau'))
&= \sum_{n \in \N_0} 
\bra \rho(-p, q, \tau) G_0, f_n \ket_{L^2(\R^D)}
\bra f_n, \rho(-p', q', \tau') G_0 \ket_{L^2(\R^D)}
\nonumber\\
&= \sum_{n \in \N_0} 
\bra \rho(p, q, \tau) f_n, G_0 \ket_{L^2(\R^D)}
\bra G_0, \rho(p', q', \tau') f_n \ket_{L^2(\R^D)}
\nonumber\\
&= \sum_{n \in \N_0} M_{f_n, G_0}(p, q, \tau)
\overline{M_{f_n, G_0}(p', q', \tau')}
\nonumber\\
&= \frac{e^{-(|x|^2+|x'|^2)/2}}{\pi^D}
\sum_{n  \in \N_0} \sB[f_n](x) \overline{\sB[f_n](x')}.
\label{eqn:Cor1}
\end{align}

The Hermite polynomials on $\R$ are defined by
\[
H_n(\zeta) := (-1)^n e^{\zeta^2}
\frac{d^n}{d \zeta^n} e^{-\zeta^2},
\quad n \in \N_0, \quad \zeta \in \R,
\]
and then the \textit{Hermite orthonormal functions} are given by
\begin{equation}
\psi_n(\zeta) := \frac{1}{\sqrt{2^n n! \sqrt{\pi}}}
e^{-\zeta^2/2} H_n(\zeta),
\quad n \in \N_0, \quad \zeta \in \R;
\label{eqn:Hermite}
\end{equation}
$\int_{\R} \psi_n(\zeta) \psi_m(\zeta) d \zeta=\delta_{nm},
\, n, m \in \N_0$.
Here we make a slight modification as
\[
\widetilde{\psi}_n(\zeta) :=
2^{1/4} \psi_n(\sqrt{2} \zeta), \quad 
n \in \N_0, \quad \zeta \in \R.
\]
$\{\widetilde{\psi}_n\}_{n \in \N_0}$ makes
a \textit{real} CONS of $L^2(\R)$, 
$\int_{\R} \widetilde{\psi}_n(\xi) \widetilde{\psi}_m(\xi) 
d \xi=\delta_{nm},
\, n, m \in \N_0$, as assumed above.

We extend the above results on $\R$ to
$\R^D$, $D \in \N$ by simply considering 
direct products;
for $n :=(n^{(1)}, \dots, n^{(D)}) \in \N_0^D$,
$\zeta:=(\zeta^{(1)}, \dots, \zeta^{(D)}) \in \R^D$,
\begin{equation}
\Psi_n(\zeta) := \prod_{\ell=1}^D 
\widetilde{\psi}_{n^{(\ell)}}(\zeta^{(\ell)}).
\label{eqn:Psi1}
\end{equation}
We can show that
\begin{align*}
\sB[\Psi_n](x)
&:= \prod_{\ell=1}^D \sB[\widetilde{\psi}_{n^{(\ell)}}](x^{(\ell)})
\nonumber\\
&= \prod_{\ell=1}^D \phi_{n^{(\ell)}}(x^{(\ell)})
=:\Phi_n(x),
\quad n \in \N_0^D, \quad x \in \C^D,
\end{align*}
where
\begin{equation}
\phi_n(x) := \frac{x^n}{\sqrt{n!}}, \quad
n \in \N_0, \quad x \in \C.
\label{eqn:phin}
\end{equation}
Hence (\ref{eqn:Cor1}) is calculated as
\begin{align*}
{\rm Cor}_{G_0}((p,q,\tau), (p', q', \tau'))
&= \frac{e^{-(|x|^2+|x'|^2)/2}}{\pi^D}
\sum_{n \in \N_0^D}
\Phi_n(x) \overline{\Phi_n(x')}
\nonumber\\
&= \frac{e^{-(|x|^2+|x'|^2)/2}}{\pi^D}
\prod_{\ell=1}^D 
\sum_{n^{(\ell)}=0}^{\infty} \frac{1}{n^{(\ell)}!}
\Big( x^{(\ell)} \overline{{x'}^{(\ell)}} \Big)^{n^{(\ell)}}
\nonumber\\
&= \sqrt{\frac{e^{-|x|^2}}{\pi^D}}
K_{\sH_D}(x, x')
\sqrt{\frac{e^{-|x'|^2}}{\pi^D}},
\quad x, x' \in \C^D.
\end{align*}
By the gauge invariance of DPP
(see Lemma \ref{thm:Gauge} below), 
this can be identified with the
correlation kernel $K_{\sH_D}(x, x')$ of
the Heisenberg family of DPPs
given by (\ref{eqn:K_HD}). 
Note that this result is obtained due to
the special choice of $g$ given by (\ref{eqn:G0}), which is
also written as
\begin{equation}
g(\zeta)=
G_0(\zeta) :=\frac{1}{\pi^{D/2}} \Psi_0(\zeta),
\quad \zeta \in \R^D.
\label{eqn:G0_2}
\end{equation}
We think that (\ref{eqn:Psi1}) represents
the eigenstate in a multivariate level
$n \in \N_0^D$ of the $D$-dimensional 
direct-product system of harmonic oscillators. 
Since $\Psi_0(\zeta)$ gives the ground state,
the above results on the choice (\ref{eqn:G0_2})
can be regarded as the 
\textit{ground-state representation}. 

\subsection{Higher-level representations}
\label{sec:higher}

Now we consider the following general choice 
of $g$,
\begin{align*}
g(\zeta) &=G_m(\zeta) :=
\frac{1}{\pi^{D/2}} \Psi_m(\zeta),
\nonumber\\
&
\quad m=(m^{(1)}, \dots, m^{(D)}) \in \N_0^D, 
\quad \zeta=(\zeta^{(1)}, \dots, \zeta^{(D)}) \in \R^D, 
\end{align*}
which corresponds to the higher-level state of
the $D$-dimensional system of harmonic oscillators. 
The matrix-element function is then given by
\[
M_{f, G_m}(p, q, \tau)
:=\bra \rho(p, q, \tau) f, G_m \ket_{L^2(\R^D)}
=e^{2 \sqrt{-1} \tau} \sB_m[f](x, \overline{x})
\frac{e^{-|x|^2/2}}{\pi^{D/2}},
\]
where
\begin{align*}
\sB_m[f](x, \overline{x}) &:= \left( \frac{2}{\pi} \right)^{D/4}
\frac{1}{\sqrt{2^m m!}}
\int_{\R^D} f(\zeta) e^{2 \zeta \cdot x- \zeta^2 -x^2/2}
\bH_m \left( \sqrt{2}
\left( \zeta-\frac{x+\overline{x}}{2} \right) \right) d \zeta
\nonumber\\
&= \frac{1}{\sqrt{m!}} e^{|x|^2} 
\frac{\partial^m}{\partial x^m} \Big[
e^{-|x|^2} \sB[f](x) \Big].
\end{align*}
Here we have used the multivariate notations,
\begin{align*}
& 2^m := \prod_{\ell=1}^D 2^{m^{(\ell)}},
\quad m! := \prod_{\ell=1}^D m^{(\ell)} !,
\quad
\frac{\partial^m}{\partial x^m}
:=\prod_{\ell=1}^D \frac{\partial^{m^{(\ell)}}}
{\partial {x^{(\ell)}}^{m^{(\ell)}}},
\quad
\bH_m(\xi):=\prod_{\ell=1}^D H_{m^{(\ell)}}(\xi^{(\ell)}).
\end{align*}
By the definition (\ref{eqn:Bargmann}), 
$\sB_0[f](x)= \sB[f](x)$. 
The integral transformation 
$\sB_m[f](x, \overline{x})$, $m \in \N_0^D$ is called the
\textit{true-$m$-Bargmann transform}
by Vasilevski \cite{Vas00},
the \textit{polyanalytic Bargmann transform} 
by Abreu and Feichtinger \cite{AF14}, and 
has been studied as the \textit{coherent state transform}
by Mouayn \cite{Mou11}.

We can read from 
Shirai's paper \cite{Shi15} that, 
for (\ref{eqn:phin}), 
\[
\frac{1}{\sqrt{m!}} e^{|x|^2}
\frac{\partial^m}{\partial x^m}
\Big[ e^{-|x|^2} \phi_n(x) \Big]
= \sqrt{ \frac{m!}{n!}} 
L_m^{(n-m)}(|x|^2) x^{n-m},
\quad n, m \in \N_0, \quad x \in \C,
\]
where $L_n^{(\alpha)}$ is defined by (\ref{eqn:Laguerre}).
Hence, we have
\[
\sB_m[\Phi_n](x, \overline{x})
= \prod_{\ell=1}^D
\sqrt{\frac{m^{(\ell)} !}{n^{(\ell)} !}} 
L_{m^{(\ell)}}^{(n^{(\ell)}-m^{(\ell)})}(|x^{(\ell)}|^2) x^{n^{(\ell)}-m^{(\ell)}}
\]
for
$m=(m^{(1)}, \dots, m^{(D)}),
n=(n^{(1)}, \dots, n^{(D)}) \in \N_0^D, \, 
x=(x^{(1)}, \dots, x^{(D)}) \in \R^D$. 
It was also shown in \cite{Shi15} that
\[
\sum_{n=0}^{\infty} \frac{1}{n!} L_m^{(n-m)}(|x|^2)
L_m^{(n-m)}(|x'|^2) (x \overline{x'})^{n-m}
=\frac{1}{m!} e^{x \cdot \overline{x'}} L_m(|x-x'|^2),
\quad m \in \N_0, \quad x, x' \in \C.
\]

For each multivariate level $m \in \N_0^D$, 
as an extension of (\ref{eqn:CorG01}), 
the correlation of two points
$(p,q,\tau)$ and $(p', q', \tau')$ in the space
of $\sH_D$ shall be given by
\[
{\rm Cor}_{G_m}((p,q,\tau), (p', q', \tau'))
:= \bra \rho(-p, q, \tau) G_m,
\rho(-p', q', \tau') G_m \ket_{L^2(\R^D)}.
\]
It is evaluated as
\begin{align*}
{\rm Cor}_{G_m}((p,q,\tau), (p', q', \tau'))
&= \sum_{n \in \N_0^D} 
M_{\Phi_n, G_m}(p, q, \tau)
\overline{M_{\Phi_n, G_m}(p', q', \tau')}
\nonumber\\
&= \frac{e^{-(|x|^2+|x'|^2)/2}}{\pi^D}
\sum_{n \in \N_0^D} 
\sB_m[\Phi_n](x, \overline{x})
\overline{\sB_m[\Phi_n](x', \overline{x'})}
\nonumber\\
&= \frac{e^{-(|x|^2+|x'|^2)/2}}{\pi^D}
e^{x \cdot \overline{x'}} 
\prod_{\ell=1}^D L_{m^{(\ell)}}(|x^{(\ell)}-{x'}^{(\ell)}|^2)
\nonumber\\
&= \sqrt{ \frac{e^{-|x|^2}}{\pi^D}}
K_{\sH_D}^{(m)}(x, x') 
\sqrt{ \frac{e^{-|x'|^2}}{\pi^D}},
\quad x, x' \in \C^D.
\end{align*}
By Lemma \ref{thm:Gauge},
this can be identified with the
correlation kernel $K_{\sH_D}^{(m)}(x, x')$ of
the extended Heisenberg family of DPPs
given by (\ref{eqn:K_HD_m}). 

\SSC
{Preliminaries and Proofs} \label{sec:preliminaries_proofs}

\subsection{Preliminaries for point processes}
\label{sec:preliminaries}

The configuration space of \textit{point process} $\Xi=\Xi(\cdot)$
is given by
\[
\Conf(S)
=\Big\{ \xi = \sum_i \delta_{x_i} : \mbox{$x_i \in S$,
$\xi(\Lambda) < \infty$ for all bounded set $\Lambda \subset S$} 
\Big\}.
\]
Let
$\cB_{\rm c}(S)$ be the set of 
all bounded measurable 
complex functions on $S$ of compact support, 
and for 
$\xi \in \Conf(S)$ and $\phi \in \cB_{\rm c}(S)$
we set 
\[
\bra \xi, \phi \ket :=\int_S \phi(x) \, \xi(dx) 
=\sum_i \phi(x_i).
\]
Random variables written in this form
are generally called
\textit{linear statistics} of a point process $\Xi$ \cite{For10}.
For a point process $\Xi$, 
if there exists a non-negative measurable function
$\rho_1$ such that
\[
\bE[\bra \Xi, \phi \ket ]
=\int_S \phi(x) \rho_1(x) \lambda(dx)
\quad \forall \phi \in \cB_{\rm c}(S).
\label{eqn:1_correlation}
\]
$\rho_1$ is called the \textit{first correlation function} of $\Xi$ 
with respect to the reference measure $\lambda(dx)$.
By definition, $\rho_1(x)$ gives the \textit{density of point}
at $x \in S$ with respect to $\lambda(dx)$.
For $n \in \N$, 
from $\xi \in \Conf(S)$ we define
\[
\xi_n := \sum_{i_1, \dots, i_n : i_j \not= i_k, j \not=k} 
\delta_{x_{i_1}} \cdots \delta_{x_{i_n}},
\]
and denote the $n$-product measure of $\lambda$ as
$\lambda^{\otimes n}(dx_1 \cdots dx_n)
:= \prod_{i=1}^n \lambda(dx_i)$. 
For a point process $\Xi$, 
if there exists a symmetric, non-negative measurable function
$\rho_n$ on $S^n$ such that
\[
\bE[\bra \Xi_n, \phi \ket ]
=\int_{S^n} \phi(x_1, \dots, x_n) 
\rho_n(x_1, \dots, x_n) \lambda^{\otimes n} (dx_1 \cdots dx_n)
\quad \forall \phi \in \cB_{\rm c}(S^n),
\]
then we say that $\rho_n$ is the 
\textit{$n$-th correlation function} of $\Xi$ 
with respect to $\lambda(dx)$.

Determinantal point process (DPP) 
is defined as follows. 
\begin{df}
A point process $\Xi$ 
on $(S, \cB_{\rm c}(S), \lambda(dx))$ is said to be a DPP with 
a measurable kernel $K : S \times S \to \C$, if 
the correlation functions 
with respect to $\lambda(dx)$ are given by
\begin{equation}
\rho_n(x_1, \dots, x_n) = \det_{1 \leq i, j \leq n}
[ K(x_i, x_j) ]
\quad \mbox{for every $n \in \N$
and any $x_1, \dots, x_n \in S$}.
\label{eqn:rho}
\end{equation}
The integral kernel $K$ is called the
\textit{correlation kernel}. 
The DPP is specified by the 
\textit{triplet $(\Xi, K, \lambda(dx))$}.
\end{df}
The following fact is well known.
\begin{lem}
\label{thm:Gauge}
Consider a non-vanishing function
$f: S \to \C$.
Even if the correlation kernel $K(x,y)$ is transformed as
\begin{equation}
K(x, y) \to
K_f(x,y):= f(x) K(x,y) \frac{1}{f(y)},
\quad x, y \in S,
\label{eqn:Gauge}
\end{equation}
all correlation functions (\ref{eqn:rho}) are
the same and hence
\[
(\Xi, K, \lambda(dx)) \law= (\Xi, K_f, \lambda(dx)).
\]
\end{lem}
The transformation (\ref{eqn:Gauge}) is
called the \textit{gauge transformation}
and the above property of DPP is
referred to \textit{gauge invariance}.
See, for instance, Lemma 3.8 in Section 3.6 of \cite{Kat16}.

If the point process is a DPP, $(\Xi, K, \lambda(dx))$, then
\begin{align*}
\bE[\bra \Xi, \phi \ket]
&= \int_{S} \phi(x) K(x, x) \lambda(dx),
\nonumber\\
\var[\bra \Xi, \phi \ket]
&= \frac{1}{2} \int_{S \times S} 
|\phi(x)-\phi(y)|^2 K(x, y) K(y,x) \lambda^{\otimes 2}(dx dy),
\quad \phi \in \cB_{\rm c}(S).
\end{align*}
In particular, when $\phi$ is the indicator function
of a bounded domain $\Lambda \subset S$; 
$\phi(x)={\bf 1}_{\Lambda}(x):=1$, if
$x \in \Lambda$, and $:=0$, otherwise, 
\begin{align*}
\bE[\bra \Xi(\Lambda) \ket]
&= \int_{\Lambda} K(x, x) \lambda(dx),
\nonumber\\
\var[\bra \Xi(\Lambda) \ket]
&= \int_{\Lambda} 
\int_{S \setminus \Lambda}
K(x, y) K(y,x) \lambda(dx) \lambda(dy).
\end{align*}

\subsection{Comments on Theorem \ref{thm:main2} }
\label{sec:MKS}

The \textit{Bessel function of the first kind} 
and the \textit{modified Bessel function of the first kind}
are defined as 
\begin{align*}
J_{\nu}(z) &:= \left( \frac{z}{2} \right)^{\nu}
\sum_{n=0}^{\infty} (-1)^n \frac{(z/2)^{2n}}{n! \Gamma(\nu+n+1)},
\nonumber\\
I_{\nu}(z) &:=\left( \frac{z}{2} \right)^{\nu}
\sum_{n=0}^{\infty} \frac{(z/2)^{2n}}{n! \Gamma(\nu+n+1)},
\quad \nu > -1, 
\quad z \in \C \setminus (-\infty, 0],
\end{align*}
respectively.
The following were proved in \cite{MKS21}.
\begin{prop}
\label{thm:H_DPP2}
For the Heisenberg family of DPPs,  
$(\Xi_{\sH_D}, K_{\sH_D}, \lambda_{\rN(0, 1; \C^D)}(dx))$ 
on $\C^D$,
$D \in \N$, the following hold,
\begin{align*}
\bE[\Xi_{\sH_D}(\B^{(2D)}_R)]
&= \frac{R^{2D}}{D!},
\\
\var[\Xi_{\sH_D}(\B^{(2D)}_R)]
&=\frac{2 R^{2D}}{(D-1)!} 
\int_0^{\infty} 
\frac{J_D(\kappa R)^2}{\kappa}
(1-e^{-\kappa^2/4}) d \kappa
\nonumber\\
&= \frac{R^{2D} e^{-2 R^2}}{D!} 
\sum_{n=0}^{D-1} 
[I_n(2R^2) +I_{n+1}(2R^2)], \quad R >0.
\end{align*}
\end{prop}

Theorem \ref{thm:main2} is concluded from the above proposition,
if we use the following
asymptotic formula of the modified Bessel functions
(see, for instance, \cite[Section 10.17]{NIST10}),
\[
I_{\nu}(x) \sim
\frac{e^x}{\sqrt{2 \pi x}}
\sum_{k=0}^{\infty} (-1)^k \frac{\alpha_k(\nu)}{k! 2^{3k}} x^{-k},
\quad \mbox{as $x \to \infty$}.
\]

\subsection{Proof of Theorem \ref{thm:polydisk1} }
\label{sec:proof}
For $R>0, n, m \in \N_0$, let
\[
p_n^{(R, m)}:= \frac{m!}{n!}  \int_0^{R^2} u^{n-m} e^{-u}
|L_m^{(n-m)}(u)|^2 du.
\]
We introduce a series of random variables 
$Y^{(R, m^{(\ell)})}_{n^{(\ell)}} \in \{0,1\}$, $m^{(\ell)} \in \N_0$, 
$n^{(\ell)} \in \N_0$, $\ell=1, \dots, D$
such that they are mutually independent and
\[
Y^{(R, m^{(\ell)})}_{n^{(\ell)}} \sim
\mu^{\rm Bernoulli}_{p^{(R,m^{(\ell)})}_{n^{(\ell)}}}, 
\]
where the right-hand side denotes
the \textit{Bernoulli measure} of probability 
$p_{n^{(\ell)}}^{(R, m^{(\ell)})}$.
By the general theory of the 
\textit{duality relations between DPPs} \cite[Theorem 2.6]{KS21},
we can prove that
\[
\Xi_{\sH_D}^{(m)}(\Delta^{(D)}_R)
\dis= \sum_{n \in \N_0^D} \prod_{\ell=1}^D Y^{(R,m^{(\ell)})}_{n^{(\ell)}}.
\] 
Then it is easy to verify that
\[
\bE[\Xi_{\sH_D}^{(m)}(\Delta^{(D)}_R)]
=\sum_{n \in \N_0^D} \prod_{\ell=1}^D p_{n^{(\ell)}}^{(R, m^{(\ell)})}
=\prod_{\ell=1}^D \sum_{k=0}^{\infty} p_k^{(R, m^{(\ell)})},
\]
and
\begin{align*}
\var[\Xi_{\sH_D}^{(m)}(\Delta^{(D)}_R)]
&= \var \left[
\sum_{n \in \N_0^D} \prod_{\ell=1}^D Y^{(R, m^{(\ell)})}_{n^{(\ell)}} \right]
=  \sum_{n \in \N_0^D} \var \left[
\prod_{\ell=1}^D Y^{(R, m^{(\ell)})}_{n^{(\ell)}} \right]
\nonumber\\
&= \sum_{n \in \N_0^D}
\left[ \prod_{\ell=1}^{D} p_{n^{(\ell)}}^{(R, m^{(\ell)})}
-\left( \prod_{\ell=1}^{D} 
p_{n^{(\ell)}}^{(R, m^{(\ell)})} \right)^2 \right]
\nonumber\\
&=\prod_{\ell=1}^D \sum_{k=0}^{\infty} p_k^{(R, m^{(\ell)})}
-\prod_{\ell=1}^D \sum_{k=0}^{\infty} (p_k^{(R, m^{(\ell)})})^2.
\end{align*}
Hence, we obtain
\[
\frac{\var[\Xi_{\sH_D}^{(m)}(\Delta^{(D)}_R)]}
{\bE[\Xi_{\sH_D}^{(m)}(\Delta^{(D)}_R)]}
=1 - 
\prod_{\ell=1}^D \frac{\sum_{k=0}^{\infty} (p_k^{(R,m^{(\ell)})})^2}
{\sum_{k=0}^{\infty} p_k^{(R, m^{(\ell)})}}, \quad m \in \N_0^D.
\]
When $D=1$, the above gives
\[
\frac{\var[\Xi_{\sH_1}^{(m)}(\B^{(2)}_R)]}
{\bE[\Xi_{\sH_1}^{(m)}(\B^{(2)}_R)]} 
\equiv
\frac{\var[\Xi_{\sH_1}^{(m)}(\Delta^{(1)}_R)]}
{\bE[\Xi_{\sH_1}^{(m)}(\Delta^{(1)}_R)]}
= 1- \frac{\sum_{k=0}^{\infty} (p_k^{(R,m)})^2}
{\sum_{k=0}^{\infty} p_k^{(R,m)}}, \quad m \in \N_0.
\]
where we have used the obvious fact that 
$\Delta_R^{(1)}=\D_R \subset \C$ and it is identified
with $\B^{(2)}_R \subset \R^2$.
Hence the first equality in (\ref{eqn:Main1}) is proved.
The second equality 
is obtained by Proposition \ref{thm:Shirai1}. 
The proof is hence complete.

\vskip 1cm
\noindent{\bf Acknowledgements} \,
The present author would like to thank 
Tomoyuki Shirai and 
Shinji Koshida for useful discussion on this subject.
This work was supported by
the Grant-in-Aid for Scientific Research (C) (No.19K03674),
(B) (No.18H01124),
(S) (No.16H06338),
(A) (No.21H04432)
of Japan Society for the Promotion of Science.
It was also supported by 
the Research Institute for Mathematical Sciences (RIMS),
a Joint Usage/Research Center located in Kyoto University.
The author thanks 
Fumio Hiroshima, Itaru Sasaki and Hayato Saigo
for organizing 
the 2021 RIMS Workshop 
`Mathematical Aspects of Quantum Fields
and Related Topics', which was held online
on December 6--8, 2021. 

\begin{small}

\end{small}
\end{document}